\documentstyle[10pt]{article}
\begin{document}

\medskip
\begin{center}
\large

{\bf Generalized Weierstrass formulae, soliton equations and
Willmore surfaces}

{\bf I. Tori of revolution and the mKDV equation}

\bigskip

\normalsize

B. G. Konopelchenko

\smallskip

{\em Dipartimento di Fisica, Universita di Lecce,
73100 Lecce, Italy,  and Budker Institute of Nuclear
Physics, Novosibirsk, 630090 Russia, e-mail:konopel@le.infn.it},

\smallskip

and

\smallskip

I. A. Taimanov

\smallskip

{\em Institute of
Mathematics, Novosibirsk, 630090 Russia, e-mail: taimanov@math.nsk.su}

\end{center}

\bigskip

\noindent{\bf Abstract}:
A new approach is proposed for study structure and properties of the total
squared mean curvature $W$ of surfaces in ${\bf R}^3$. It is based on the
generalized Weierstrass formulae for inducing surfaces. The quantity
$W$ (Willmore functional or extrinsic Polyakov action) is shown to be
invariant under the modified Novikov--Veselov hierarchy of integrable
flows.

It is shown that extremals of $W$ (Willmore surfaces) obey the
two--di\-men\-si\-o\-nal Schrodinger equation and possible relations between
Willmore surfaces and the Novikov--Veselov hierarchy
of integrable equations are discussed.

The $1+1$--dimensional case and, in particular, Willmore tori of revolution,
are studied in details. The Willmore conjecture  is proved for the
mKDV--invariant Willmore tori.

\vskip1.5cm

{\bf 1. Introduction}

\medskip

Surfaces and their dynamics play an important role in many interesting
phenomena both in classical and quantum physics (see, e.g.,
\cite{BCC,Pelce,NPW,GPW}). In the string theory and $2D$--gravity based
on the Polyakov integral over surfaces (\cite{Polyakov}) contributions from
certain special classes of surfaces are of crucial importance.

In mathematics the foundations of
differential geometry of surfaces have been basically
completed almost one century ago (see, e.g., \cite{Darboux,Bianchi}).
Nevertheless, number of problems concerning the structure and properties
of special classes of surfaces are still remain open (see, e.g.,
\cite{Yau,Willmore1982,Willmore1993}).

The total squared mean curvature
$$
W = \int H^2 d\mu  \eqno({1.1})
$$
is one of the most
important characteristics of surfaces in the three-dimensio\-nal
Euclidean space ${\bf R}^3$. Here $H$ is a mean curvature and $d\mu$ is a
Liouville measure with respect to metric induced by immersion. In string
theory and $2D$--gravity the functional $W$ is known as the Polyakov's
extrinsic action. The properties of the Polyakov's action were the subject
of study in a number of papers (see the review \cite{GM}).

In the differential geometry the functional $W$ has been investigated within
the study of so-called Willmore surfaces
(\cite{Willmore1982,Willmore1993,ST,Willmore,LY,Kusner,FPPS,HJP,BB,Babich}).
These surfaces provide extremum to $W$ and obey the corresponding
Euler--Lagrange equation
$$
\Delta H + 2H(H^2-K) = 0 \eqno{(1.2)}
$$
where $K$ is the Gaussian curvature and $\Delta$ is the Laplace--Beltrami
operator. A study of extremals of the functional $W$ is also the great
importance for $2D$-gravity. In quantum theory they correspond to so-called
zero modes.

In the present paper we propose a new approach to study structure and
properties of the functional $W$ (1.1) and corresponding extremals of (1.2).
Our approach is based on the generalization of the Weierstrass formulae
proposed in \cite{Kon}. This generalization allows one to induce generic
surfaces in ${\bf R}^3$ via solutions of the system of two linear equations.
We show that the generalized Weierstrass inducing from \cite{Kon} is
equivalent to earlier extension of Weierstrass formulae given in \cite{Ken}
(see also \cite{HO,CK}). An advantage of new formulation is that it allows
to describe and construct integrable deformations of induced surfaces via
the modified Novikov--Veselov (mNV) hierarchy flows (\cite{Bogdan}, for
Veselov--Novikov equation see \cite{VN}, and also \cite{DKN}).
It is shown that the functional $W$ is invariant under mNV deformations.

In this paper (part I) we restrict ourselves to the $1+1$-dimensional limit
of our scheme. We concentrate on the mKDV--invariant Willmore tori of
revolution. Properties of such tori are studied. For mKDV--invariant
Willmore tori of revolution we prove the Willmore conjecture that their
total squared mean curvature $W \geq 2\pi^2$ and, moreover, $W=2\pi^2$
only for the Clifford torus.

The relation of this problem to soliton theory was noticed about
ten years ago and the finite--gap integration was used for
constructing Willmore tori (see \cite{FPPS,BB,Babich}).
But the formulae for finite-gap solutions is rather unefficient.
We will discuss only Willmore tori of revolution
postponing the consideration of general case and will not give rigorous
mathematical proofs (but notice that that is possible as one will see).
Our main goal is to shed light on the one-dimensional case for using this
approach to the main problem.

The paper is organized as follows. In section 2 the generalized Weierstrass
inducing is discussed. The equivalence of equation (1.2) to the
2D-Schrodinger equation is proved in section 3. The $1+1$--dimensional case
is considered in section 4. The construction
of the Clifford torus via the generalized Weierstrass inducing is given in
section 5. In section 6 we prove the Willmore conjecture for the
mKDV--invariant Willmore tori of revolution.

\newpage

{\bf 2. Generalized Weierstrass inducing}

\medskip

The generalization of the Weierstrass inducing, proposed in \cite{Kon},
starts with the linear system
$$
\psi_{1z}=p\psi_2,
$$
$$
\psi_{2\bar z}=-p\psi_1 \eqno{(2.1)}
$$
where $p(z,\bar z)$ is a real-valued function, $\psi_1$ and $\psi_2$
are complex-valued functions of the complex variable $z$, and bar
denotes the complex conjugation. Then one can define three functions
$(X^1(z,\bar z),X^2(z,\bar z),X^3(z,\bar z))$ as follows
$$
X^1+iX^2=2i\int_{\Gamma}({\bar\psi}^2_1 dz' - {\bar\psi}^2_2 d{\bar z}'),
$$
$$
X^1-iX^2=2i\int_{\Gamma}(\psi^2_2 dz' - \psi^2_1 d{\bar z}'),
\eqno{(2.2)}
$$
$$
X^3=-2\int_{\Gamma}(\psi_2{\bar\psi}_1 dz' + \psi_1{\bar \psi}_2 d{\bar z}')
$$
where $\Gamma$ is an arbitrary curve and treat $X^1,X^2$, and $X^3$ as the
coordinates of surface in ${\bf R}^3$ while $z,{\bar z}$ are local
coordinates on the surface.

An arbitrary surface in n ${\bf R}^3$ with non-vanishing mean curvature
can be represented by using this method (that was proved for Kenmotsu
representation (\cite{Ken,HO}) and valid for (2.2) since this representation
is equivalent to Kenmotsu's one as it will be shown below).

The metric which is induced by this immersion is equal to (\cite{Kon})
$$
4u(z,{\bar z})^2 dz d{\bar z}   \eqno{(2.3)}
$$
where $u(z,{\bar z})=|\psi_1|^2+|\psi_2|^2$, and the Gaussian curvature
$K$ and the mean curvature $H$ have the form
$$
K=-\frac{1}{4u^2}\Delta \log u, \ \  H=\frac{p}{u}.   \eqno{(2.4)}
$$

It is easy to notice that formulae (2.1) and (2.2) are reduced to those of
Weierstrass in the case $p=0$, i.e., for minimal surfaces.

Exact explicit solutions to (2.1) provide us surfaces in ${\bf R}^3$
via (2.2) by quadratures. A family of surfaces parametrized by an arbitrary
number of holomorphic functions has been constructed in \cite{Kon}.

Another extension of the Weierstrass inducing to nonminimal surfaces has
been proposed earlier  by Kenmotsu (\cite{Ken}, see also \cite{HO}). In
this approach a surface in ${\bf R}^3$ is defined :
$$
X^i = \int \omega_i \eqno{(2.5)}
$$
where
$$
\omega_1 = Re \{ \phi (1-f^2) dz\},
$$
$$
\omega_2 = Re \{i\phi(1+f^2) dz\}, \eqno{(2.6)}
$$
$$
\omega_3 = Re \{ 2\phi f dz \}
$$
and functions $f$ and $\phi$ obey the relation
$$
(\log \phi)_{\bar z} = - \frac{2{\bar f}f_{\bar z}}{1+|f|^2}. \eqno{(2.7)}
$$
The mean curvature has the form
$$
H = -\frac{{\bar f}_{\bar z}}{\phi(1+|f|^2)}.
$$

The generalized Weierstrass inducing (2.1) --- (2.2) and the Kenmotsu
inducing (2.5) --- (2.7) are equivalent to each other. The relation
between $\psi_1,\psi_2,p$ and $f,\phi$ are of the form (see also \cite{CK})
$$
f= i\frac{{\bar \psi}_1}{\psi_2}, \phi = i\psi^2_2 \eqno{(2.8)}
$$
and
$$
p= -\frac{({\bar f}{\bar \phi})_z}{\sqrt{\phi{\bar \phi}}} =
\frac{\phi f_{\bar z}}{\sqrt{\phi{\bar \phi}}(1+|f|^2)}. \eqno{(2.9)}
$$

One of advantages of the representation (2.1) -- (2.2) is that, due to
the linear character of (2.1), it allows to describe integrable
deformations of induced surfaces in a simple manner. For that purpose
one uses the basic idea of the inverse spectral transform method and
looks for the deformation of all quantities in (2.1) in time $t$ via
the equation
$$
\left(\begin{array}{c} \psi_1 \\\psi_2 \end{array} \right)_t =
A \left(\begin{array}{c} \psi_1 \\ \psi_2 \end{array}
\right) \eqno{(2.10)}
$$
where $A$ is a differential $2 \times 2$--matrix operator (\cite{Kon}).
The requirement of compatibility of (2.1) and (2.10) first quarantees the
preservation in time $t$ of the generalized Weierstrass inducing (2.2)
and, second, leads to nonlinear differential equations for $p$.

The simplest nontrivial example corresponds to nonlinear problem (2.10)
of the form
$$
(\partial_t + \partial^3_z + \partial^3_{\bar z} +
3\left(\begin{array}{cc} 0 & p_z \\ 0 & \omega \end{array}\right)\partial_z
+ 3\left(\begin{array}{cc} {\bar \omega} & 0 \\ p_{\bar z} & 0
\end{array}\right)\partial_{\bar z} +
$$
$$
\frac{3}{2}\left(\begin{array}{cc} {\bar \omega}_{\bar z} & 2p\omega \\
-2p{\bar \omega} & \omega_z \end{array} \right)) \psi = 0.
\eqno{(2.11)}
$$

The associated nonlinear integrable equation for $p$ is
$$
p_t + p_{zzz} + p_{{\bar z}{\bar z}{\bar z}} +
3p_z\omega + 3p_{\bar z}{\bar \omega} +
\frac{3}{2}p{\bar \omega}_{\bar z} + \frac{3}{2}p\omega_z = 0,
\eqno{(2.12)}
$$
$$
\omega_{\bar z} = (p^2)_z.
$$
Equation (2.12) is known as the modified Novikov--Veselov equation
(mNV) (\cite{Bogdan}). The mNV equation is integrable by the Inverse
Scattering Transform (IST) method. The hierarchy of integrable PDEs
associated with the system (2.1) arises as the compatibility condition
of (2.1) with the system of the form (2.10) where
$$
A = - (\partial^{2k+1}_z + \partial^{2k+1}_{\bar z}
+ \sum_{m=1}^{2k-1} (u_m\partial^m_z + {\bar u}_m \partial^m_{\bar z})).
\eqno{(2.13)}
$$
All equations of mNV hierarchy commute to each other and are
integrable by the IST method.

Thus the integrable dynamics of surfaces is induced by the mNV hierarchy
via (2.2). This integrable dynamics of surfaces inherits all properties
of the mNV hierarchy. Note that the minimal surfaces ($p=0$) are
invariant under such dynamics.

Within the generalized Weierstrass inducing the functional $W$ has a very
simple form. Indeed, by using of (2.3) and (2.4), one gets
$$
W= 4 \int p^2 dz d{\bar z}. \eqno{(2.14)}
$$
Thus, using the exact solutions of system (2.1), one gets $W$.

An important properties of $W$ is that it is invariant under the
mNV deformations. Indeed, rewriting the mNV equation (2.12) in the form
$$
(p^2)_t + (2p_{zz} - p^2_z + 3p^2\omega)_z +
(2p_{{\bar z}{\bar z}} - p^2_{\bar z} + 3p^2{\bar \omega})_{\bar z} = 0
$$
one easily concludes that
$$
W_t = 4 \int (p^2)_t dz d{\bar z} = 0    \eqno{(2.15)}
$$
for periodic surfaces ($p$ is a periodic function) or for $p$ decreasing as
$|z| \rightarrow \infty$. The quantity $\int p^2 dz d{\bar z}$ is the
integral of motion for the whole mNV hierarchy (we will give the complete
prove for $1+1$--case below).

Thus the infinite--parametric family of mNV deformed surfaces have
the same value of $W$.

\vskip1.5cm

{\bf 3. Willmore surfaces, two--dimensional Schrodinger equation and
NV hierarchy}

\medskip

We consider now the Euler--Lagrange equation  (1.2) which defines the
Willmore surfaces. We will consider a surface in the system of local
coordinates formed by minimal lines. The metric in these coordinates has
the form (2.3) where $u$ is some real--valued function.
Then let us  introduce the real--valued function $p$ via $H=\frac{p}{u}$.

Since in these coordinates
$$
H=\frac{H_{z{\bar z}}}{u^2}, K = -\frac{(\log u)_{z{\bar z}}}{u^2}
$$
equation (1.2) becomes
$$
u_{z{\bar z}} - (\log p)_{\bar z}u_z - (\log p)_zu_{\bar z} +
\frac{p_{z{\bar z}}+2p^3}{p}u =0. \eqno{(3.1)}
$$
Then we introduce new variable $\xi$ via $u=p\xi$. In terms of $\xi$
equation (3.1) is of the form
$$
\xi_{z{\bar z}} + 2((\log p)_{z{\bar z}} + p^2) \xi =0. \eqno{(3.2)}
$$
Noting that $\xi = \frac{u}{p} = \frac{1}{H}$, we finally conclude that
the Euler--Lagrange equation (1.2) is equivalent to the two--dimensional
Schrodinger equation
$$
(\frac{1}{H})_{z{\bar z}} + V\frac{1}{H} = 0 \eqno{(3.3)}
$$
with the potential
$$
V = 2( \log(uH))_{z{\bar z}} + 2(uH)^2. \eqno{(3.4)}
$$
In terms of the variable $\log p = \log (uH) = \phi$ one has
$$
V = 2(\phi_{z{\bar z}} + \exp 2\phi).   \eqno{(3.5)}
$$

Thus for induced Willmore surfaces the linear system
(3.3) holds.

Note that for periodic $p$ and $p$ decaying as $|z| \rightarrow \infty$,
one has
$$
W= 4 \int p^2 dz d{\bar z} = 2\int V dz d{\bar z}.  \eqno{(3.6)}
$$

Following the basic idea of the Inverse Scattering Transform method
(see \cite{Nov,AS}), we define deformations via the compatible system of
equations
$$
\xi_{z{\bar z}} + V\xi=0,
$$
$$
\xi_t+A(\partial_z,\partial_{\bar z},V)\xi=0 \eqno{(3.7)}
$$
where $A$ is a differential operator and $t$ is a deformation
parameter. The compatibility condition (3.7) is equivalent to nonlinear PDE
for $V$. The simplest nonlinear PDE of this type is the Novikov--Veselov
equation (\cite{VN})
$$
V_t + V_{zzz} + V_{{\bar z}{\bar z}{\bar z}} +
3(V(\partial_{\bar z}^{-1}V_z))_z +
3(V(\partial_z^{-1}V_{\bar z}))_{\bar z} = 0.  \eqno{(3.8)}
$$
This equation is the first one from the infinite hierarchy of
Novikov--Veselov equations which are integrable by the IST method.
It is easy to show that for the Novikov--Veselov equation (3.8) as for
all equations from this hierarchy the identity
$$
\frac{\partial}{\partial t}\int V dz d{\bar z} = 0.  \eqno{(3.9)}
$$
holds.

Notice that in one--dimensional limit the equation (3.3) reduces to the
linear problem for one--dimensional Schrodinger equation and the
Novikov--Veselov equation (3.8) reduces to the Korteweg--de Vries equation
$$
V_t + \frac{1}{4}V_{xxx} + 12VV_x = 0. \eqno{(3.10)}
$$

Possible geometrical meaning of such deformations
will be considered elsewhere.

\vskip1.5cm

{\bf 4. Induced tori of revolution and the mKDV equations}

\medskip

In the rest of the paper we will consider surfaces of revolution
that correspond to the case when
$p(z,{\bar z})=p(x)$ where $z=x+iy$ and $x$ and $y$ are real-valued
variables, and functions $\psi_1$ and  $\psi_2$ which define the
generalized Weierstrass inducing have the following form
$$
\psi_1=r(x) \exp{\frac{iy}{2}} , \ \   \psi_2=s(x) \exp{\frac{iy}{2}}
\eqno{(4.1)}
$$
where $r(x)$ and $s(x)$ are real-valued functions.

For surfaces of revolution system (2.1) takes the form
$$
r_x=-\frac{1}{2}r +2ps,
$$
$$
s_x=\frac{1}{2}s-2pr  \eqno{(4.2)}
$$

These equation means that the vector function $(r,s)$
belongs to the kernel of the following operator
$$
L = \partial_x - \frac{1}{2}
\left(\begin{array}{cc} \lambda & v \\ -v & -\lambda \end{array}\right)
\eqno{(4.3)}
$$
for
$$
\lambda=-1, v=4p.
$$

Since (4.3) we have
$$
(rs)_x=-\frac{v}{2}(r^2-s^2), (r^2+s^2)_x=\lambda(r^2-s^2),
$$
$$
(r^2-s^2)_x = \lambda(r^2+s^2)+2vrs.          \eqno{(4.4)}
$$
It follows directly from (4.4) that functions $u$ and $p$ are related
via the following identity
$$
v=\frac{\lambda^2u-u_{xx}}{\sqrt{\lambda^2u^2-u_x^2}}
\eqno{(4.5)}
$$
which for $\lambda=-1, p=4v$ has the form
$$
p=\frac{u-u_{xx}}{4\sqrt{u^2-u^2_x}}\eqno{(4.6)}
$$

We call the induced surface of revolution one-periodic  if
the functions $r(x)$ and $s(x)$ are periodic with the same period $T$.
In this case it follows from (2.2) that functions $X^1$ and $X^2$ are
double--periodic but the formula for $X^3$ has the form
$$
X^3 = -4\int r(x)s(x) dx  \eqno{(4.7)}.
$$
Generally formulae (2.2) induce in this case a cylinder
of revolution.
We also mention that if
$$
\int_0^T r(x)s(x) dx = 0 \eqno{(4.8)}
$$
then in virtue of (4.7) this one-periodic cylinder converts into a torus
of revolution.
Via (2.3) and (2.4) the total squared mean
curvature of such torus $M$ is equal to
$$
W(M)=\int_M H^2 d\mu = \int_0^T dx \int_0^{2\pi} dy \frac{p^2}{u^2} 4u^2  =
8\pi \int_0^T p^2 dx           \eqno{(4.9)}
$$
where $d\mu = 4u^2 dx dy $ is the natural measure on $M$ which is determined
the induced metric.

Now let us introduce mKdV hierarchy. We will follow formulae from \cite{S}.

Let us consider the following matrix operator:
$$
K_{2n+1}=\frac{1}{2}\left(\begin{array}{cc}
A & B \\ C & -A \end{array}\right) \eqno{(4.10)}
$$
where matrix elements are determined by following formulae
$$
A=\sum_{k=0}^{n} A_{2k+1}^{(n)}\lambda^{2k+1},
$$
$$
B+C=2 \sum_{k=0}^{n-1}S_{2k+1}^{(n)}\lambda^{2k+1},
B-C=2\sum_{k=0}^{n} T_{2k}^{(n)} \lambda^{2k},
$$
$$
A_{2k+1}^{(n)} = \partial^{-1}_x (vD^{n-k-1}v_x), A_{2n+1}^{(n)} = 1,
\eqno{(4.11)}
$$
$$
S_{2k+1}^{(n)}=D^{n-k-1}v_x, \ \ T_{2k}^{(n)} = \partial^{-1}_x D^{n-k}v_x,
$$
$$
D= \partial^2_x + v^2 + v_x \partial^{-1}_x v.
$$

The Lax equations
$$
[\partial_{t_{2n+1}} + K_{2n+1}, L] = 0 \eqno{(4.12)}
$$
form hierarchy of nonlinear equations known as modified Korteweg--de Vries
(mKdV) hierarchy. These equations are equivalent to the following
$$
v_{t_{2n+1}} = D^n v_x  \eqno{(4.13)}
$$
and the first nontrivial equation among them is the famous mKdV--equation
$$
v_t=\frac{3}{2}v^2v_{x}+v_{xxx}. \eqno{(4.14)}
$$
One can see that (4.14) is one--dimensional limit of (2.12) (for $p=
\frac{v(x)}{4})$ or (3.8).

Deformation (4.13) induces deformation of $r(x,t_1,\dots)$ and
$s(x,t_1,\dots)$
$$
\eta_{t_{2n+1}} = K_{2n+1} \eta \eqno{(4.15)}
$$
where $\eta^T = (r,s)$. Let us show that condition (4.8) is preserved by
all these flows.

Let us compute
$$
I_n = \frac{\partial}{\partial t_{2n+1}} \int^T_0 rs dx =
\int^T_0 (r_{t_{2n+1}}s+rs_{t_{2n+1}}) dx =
$$
$$
\int^T_0 \left(\sum_0^{n-1} S^{(n)}_{2k+1} \lambda^{2k+1} (r^2+s^2) +
\sum_0^n T^{(n)}_{2k} \lambda^{2k} (s^2-r^2)\right) dx. \eqno{(4.16)}
$$

Integrating (4.16) by parts with using of (4.4) we convert it as follows
$$
I_n= \int^T_0 \left(\sum_0^{n-1} \lambda^{2k+1} D^{n-k-1}v_x +
\sum_0^n \lambda^{2k-1} D^{n-k}v_x \right) (r^2+s^2) dx.
$$

Thus if we will prove that
$J_k=\int^T_0 (D^kv_x) (r^2+s^2) dx =0$
for any $k \geq 0$ then we will prove that $I_n=0$ for any $n$.

It is easy to notice that it follows from (4.4) that
$$
J_0=\int^T_0 v_x (r^2+s^2) dx = - \int_0^T v\lambda(r^2-s^2) dx = 2\lambda
\int^T_0 (rs)_x dx =0.
$$
Let us compute $J_k$ (we will omit limits of integration):
$$
J_k = \int (D^k v_x) (r^2+s^2) dx =
$$
$$
\int(\partial^2_x + v^2 +v_x\partial_x^{-1}v)
(D^{k-1}v_x) (r^2+s^2) dx =
F_1+F_2+F_3
$$
where in virtue of (4.4)
$$
F_1= \int (D^{k-1}v_x) (r^2+s^2)_{xx} dx =
$$
$$
\lambda^2\int (D^{k-1}v_x)
(r^2+s^2) dx + 2\lambda\int vrs (D^{k-1}v-x) dx,
$$
$$
F_2= \int v^2 (D^{k-1}v_x) (r^2+s^2) dx,
$$
$$
F_3= \int v_x\partial^{-1}_x(vD^{k-1}v_x)(r^ 2+s^2)dx =
\int \partial^{-1}_x(vD^{k-1}v_x) (r^2+s^2) dv=
$$
$$
- \int v^2 (D^{k-1}v_x) (r^2+s^2) dx - \lambda\int v(r^2-s^2)
\partial^{-1}_x(vD^{k-1}v-x) dx =
$$
$$
-\int v^2 (D^{k-1}v_x) (r^2+s^2) dx -2\lambda\int vrs (D^{k-1}v-x)dx.
$$
Combining all formulae for $F_1,F_2$ and $F_3$ we derive that
$$
J_k= \lambda^2 J_{k-1}
$$
and since $J_0=0$ we obtain that $J_k=0$ for all $k$.
As we mentioned above it immediately follows now that
$$
I_n =0, \ \ \ n \geq 0.
$$
Thus we prove that

{\it  equations of the mKdV hierarchy (4.13) deform tori, induced
via (2.2) and (4.1), into tori.}

Let us now show that the total squared
mean curvature $W(M)$ is also preserved
by these flows. Firstly mention that operator
$$
D^+= \partial^2_x + v^2 - v\partial^{-1}v_x
$$
is formally coadjoint to operator $D$, i.e.
$$
\int f\cdot Dg dx = \int D^+f \cdot g dx.
$$
One can see that
$$
\partial_x D^+ = D \partial_x.
$$.

Since $v=4p$, the identity
$$
W(M)=\frac{\pi}{2} \int_0^T v^2 dx \eqno{(4.17)}
$$
holds. The derivative of $W(M)$ with respect to $t_{2k+1}$
is equal to
$$
\frac{2}{\pi}W(M)_{t_{2k+1}}=2\int v (D^k v_x) dx=\int v_x ((D^+)^k v) dx =
$$
$$
(v (D^+)^k v)|_0^T - \int v \partial_x ((D^+)^k v) dx =
$$
$$
(v (D^+)^k v)|_0^T - \int v (D^k v) dx.
$$
It is only enough to mention now that the function $v (D^+)^k v$ is
periodic and thus we derive that
$$
\int  v (D^k v) dx = \frac{1}{2} (v (D^+)^k v) |_0^T =0.
$$

Thus we conclude that

{\it the total  squared mean curvature $W$ is the first integral of
all flows of the mKdV--hierarchy :
$W_{t_{2k+1}} \equiv 0$ for $k > 0$.}

Now we proceed to Willmore tori.
The Euler--Lagrange equation for the Willmore functional
in terms of these functions has the following form
$$
\frac{1}{u^2}\{\frac{1}{4}(\frac{p}{u})_{xx}+\frac{2p^3}{u} +
\frac{p}{2u}(\log{u})_{xx}\}=0.\eqno{(4.18)}
$$
Multiplying (4.18) by $4u^4$ we obtain
$$
p_{xx}u + pu_{xx} - 2p_xu_x + 8up^3 = 0.\eqno{(4.19)}
$$
Notice that (4.19) follows from (3.1) as its one--dimensional reduction.
Differentiating (4.19) by $x$ we obtain
$$
u_{xxx}p-u_{xx}p_x+u_x(-p_{xx}+8p^3) + u(24p^2p_x + p_{xxx})=0.\eqno{(4.20)}
$$

But also we can differentiate (4.6) by $x$
(for $v=4p, \lambda=-1$ ) and conclude that
$$
u_{xxx}p-u_{xx}p_x - u_x(1-16p^2)p + up_x = 0. \eqno{(4.21)}
$$
It follows from (4.20) and (4.21) that for Willmore
tori of revolution the important equation
$$
(8p^3+p_{xx} - p) u_x + (p_x -24p^2p_x - p_{xxx}) u =0.\eqno{(4.22)}
$$
holds.

It follows from (4.22) that
$$
8p^3 + p_{xx} - p = a \cdot u , a=const. \eqno{(4.23)}
$$

Thus we conclude that

{\it for every induced Willmore surface of revolution equation
(4.23) holds. }

\vskip1.5cm

{\bf 5. Clifford Torus}

\medskip

Here we demonstrate how the Clifford torus can be obtained by inducing
(2.2),(4.1).

Let $S^3$ be a unit sphere in four-dimensional Euclidean space
${\bf R}^4$. The Clifford torus in ${\bf R}^4$ is an image of the immersion
$$
{\bf R}^2 \rightarrow S^4 :
(u,v) \rightarrow (\frac{\cos{u}}{\sqrt{2}}, \frac{\sin{u}}{\sqrt{2}},
\frac{\cos{v}}{\sqrt{2}}, \frac{\sin{v}}{\sqrt{2}}). \eqno{(5.1)}
$$
It is easy to notice that this immersion is double periodic and its image
will be an embedded torus.

Let us consider the stereographic projection of $S^4$ onto the plane
$x_4=-1$ from the pole $(0,0,0,1)$:
$$
(x_1,x_2,x_3,x_4) \rightarrow (\frac{-2x_1}{x_4-1},\frac{-2x_2}{x_4-1},
\frac{-2x_3}{x_4-1},-1).
\eqno{(5.2)}
$$
We will call the image of the Clifford torus with respect to this projection
Clifford (in ${\bf R}^3$ --- !) again:
$$
(u,v) \rightarrow (\frac{2\cos{u}}{D}, \frac{2\sin{u}}{D},
\frac{2\cos{v}}{D})
\eqno{(5.3)}
$$
where $D=\sqrt{2} - \sin{v}$.
One can compute that the first fundamental form is equal to
$$
\frac{4}{D^2}( du^2 + dv^2 ) \eqno{(5.4)}
$$
and the second fundamental form is equal to
$$
\frac{2(\sqrt{2}\sin{v} - 1)}{D^2} du^2 + \frac{2}{D^2} dv^2.
\eqno{(5.5)}
$$
It follows now that the Gaussian curvature $K$ and the mean curvature $H$
are given by the following formulae
$$
K= \frac{\sqrt{2}\sin{v} - 1}{4}, \ \ H=\frac{\sin{v}}{2\sqrt{2}}.
\eqno{(5.6)}
$$

Let us now put
$$
p(x) = \frac{\sin{x}}{2\sqrt{2}(\sqrt{2}-\sin{x})}
\eqno{(5.7)}
$$
and
$$
u(x)= \frac{1}{\sqrt{2}-\sin{x}}.
\eqno{(5.8)}
$$
It is easy to check now by using of direct computations that functions
$r(x)$ and $s(x)$ such that
$$
r^2=\frac{u-u_x}{2}, s^2=\frac{u+u_x}{2},
rs=\frac{\sqrt{2}\sin{x}-1}{2(\sqrt{2}-\sin{x})^2}
\eqno{(5.9)}
$$
satisfy system (4.2).

We can also obtain functional equation for $p(x)$. It follows from (5.7)
that
$$
\sin x = \frac{8p}{1+2\sqrt{2}p}. \eqno{(5.10)}
$$
Differentiating (5.10) by $x$ we obtain expression for $\cos x$
in terms of $p$ and $p_x$:
$$
\cos x = \frac{8p_x}{(1+2\sqrt{2}p)^2}.  \eqno{(5.11)}
$$
Substituting (5.10) and (5.11) into the trivial identity
$$
\sin^2 x + \cos^2 x =1
$$
we obtain
$$
p_x^2= -4p^4+2p^2+\frac{p}{\sqrt{2}}+\frac{1}{16} \eqno{(5.12)}
$$

\vskip1.5cm

{\bf 6. mKDV--stationary Willmore tori of revolution}

\medskip

We will restrict ourselves to the Willmore tori generated by stationary
solutions of the mKDV equation (4.14) (for $v=4p$).
Such solutions obey the equation
$$
c_2p_x + p_{xxx} + 24p^2p_x =0 \eqno{(6.1)}
$$
where $c_2$ is an arbitrary constant. The relation (6.1) implies that
$$
p_{xx} + 8p^3 - c_2p - \frac{c_1}{2} = 0 \eqno{(6.2)}
$$
or finally
$$
p_x^2= -4p^4 + c_2p^2 +c_1p + c_0   \eqno{(6.3)}
$$
where $c_0$ and $c_1$ are also arbitrary constants.

It is easy to notice that this function $p(x)$ is elliptic.

Substituting (6.2) into (4.23) we obtain
$$
a \cdot u = (c_2-1)p + \frac{c_1}{2}.       \eqno{(6.4)}
$$
Thus we arrive to two possibilities

1) $ a \neq 0  $ ;

2) $a=0$ and in this case $c_2=1, c_1=0$.

In the first case substituting  (6.4) into (4.6) we obtain that
$$
c_2=2, c_0 = \frac{4c_1^2-1}{16}. \eqno{(6.5)}
$$
But one can also substitute (6.4) into (4.19) and obtain that
$$
c_1=c_2d, c_0=\frac{c_1d}{4}  \eqno{(6.6)}
$$
where $d=\frac{c_1}{2(c_2-1)}$.
Combining (6.5) and (6.6) together we obtain that
$$
c_2=0, c_1^2=\frac{1}{2}, c_0=\frac{1}{16}.
$$
But as one can see (5.12) that we obtain Clifford torus (note that
functions $p(x)$ and $-p(x)$ induce the same surface). Thus

{\it $ a \neq 0$ in (6.4) then the only Willmore surface of revolution
for which (6.4) holds is the Clifford torus.}

It is left now to consider the case when
$$
p_x^2=-4p^4+p^2+\alpha \eqno{(6.7)}
$$
where $\alpha$ is a real parameter.
This family may contain non-trivial periodic potentials of
two different types:
1) for $\frac{-1}{16} < \alpha < 0$ potential $p(x)$ is positive or negative
and
$\frac{1-\sqrt{1+16\alpha}}{8} \leq p^2 \leq\frac{1+\sqrt{1+16\alpha}}{8}$ ;
2) potential $p$ varies from $-C$ till $C$ where
$C=\sqrt{1+\frac{1+\sqrt{1+16\alpha}}{8}}$ for $\alpha > 0$.
We will show that by different reasons every of these families  does not
contain Willmore torus of revolution with $W \leq 2\pi^2 $.

Let us assume that potential $p(x)$ which satisfies (6.7) induces a
one-periodic cylinder of revolution.

It follows from (4.6) that identity (4.7) is equivalent to the
following
$$
\delta_0 = \int_0^T \frac{u-u_{xx}}{p} dx = 0. \eqno{ (6.8)}
$$
Integrating (6.8) by parts two times we obtain that
$$
\delta_0 = \int_0^T u(\frac{1}{p}+\frac{p_{xx}}{p^2}-\frac{2p_x^2}{p^3}) dx.
\eqno{(6.9)}
$$

But in virtue of (6.7),
$$
\delta_0 = -2\alpha \int_0^T \frac{u}{p^3} dx. \eqno{(6.10)}
$$
Since $u=r^2+s^2>0$ then for $\alpha <0$ integrand in the right-hand
side of (6.10) does not change its sign and we conclude that
$\delta_0 \neq 0$ and thus potentials (6.7) with $\alpha<0$ do not
induce tori of revolution.

We will not obtain an analogous statement for potentials (6.7) with
$\alpha >0$ but we will usually show that such potential induces torus
then inequality
$$
W= 8\pi \int_0^T p^2(x) dx > 2\pi^2. \eqno{(6.11)}
$$
holds.

Let us assume that
potential $p(x)$ which satisfies (6.7) with $\alpha>0$ induces
a torus via formulae (2.2) and (4.1).
Let $T$ be a period of functions $p(x), r(x)$, and $s(x)$. The total
squared mean curvature is equal to
$$
W = 8\pi\int_0^T p^2 dx =
32\pi\int_0^{\max p} \frac{p^2 dp}{\sqrt{-4p^4+p^2+\alpha}}.
\eqno{(6.12)}
$$
Let us put
$$
\beta=\sqrt{1+16\alpha}, C=\frac{1+\beta}{8}, k^2=\frac{1+\beta}{2\beta}.
$$
It is easy to notice that $\max p = -\min p =\sqrt{C}$ and substituting
$v=C-p^2$ into (6.12) we obtain $W=16\pi I$ where
$$
I= \int_0^C \frac{\sqrt{C-v} dv}{\sqrt{v(\beta-4v)}}. \eqno{(6.13)}
$$
But the last formula can be expressed in terms of classic elliptic integrals
(see \cite{GR}, formula 3.141.8):
$$
I=\sqrt{\beta}E(k)-\frac{\beta-1}{2\sqrt{\beta}}F(k) \eqno{(6.14)}
$$
where $E(k) = E(\frac{\pi}{2},k), F(k)=F(\frac{\pi}{2},k)$, and
$$
E(\phi,k)=\int_0^{\phi}\sqrt{1-k^2 \sin^2(\tau)} d\tau,
$$
$$
F(\phi,k)=\int_0^{\phi}\frac{d\tau}{\sqrt{1-k^2\sin^2(\tau)}}.
$$
The right-hand side of (6.14) is equal to
$f(k)=\frac{E(k)-(1-k^2)F(k)}{\sqrt{2k^2-1}}$ and it is easy to see that
this function is well defined for $k \in (1/\sqrt{2},1]$ and continuous
on this set, $f(1)=1$ and $f(k) \rightarrow +\infty$ as $k \rightarrow
1/\sqrt{2}$, and the inequality $W>2\pi^2$ is equivalent to the following
$$
f(k) > \frac{\pi}{4},  \frac{1}{\sqrt{2}} < k < 1. \eqno{(6.15)}
$$

Let us find the minimum of this function. Since $f(k)$ is smooth, we ought
firstly to find points where $f'(k)=0$.
But two following important identities holds (\cite{GR}, formulae
8.123.2 and 8.124.1):
$$
E(k)-(1-k^2)F(k) = k(1-k^2) \frac{dF}{dk},   \eqno{(6.16)}
$$
$$
\frac{d}{dk} \{k(1-k^2)\frac{dF}{dk}\} = kF.  \eqno{(6.17)}
$$
In virtue of (6.16) and (6.17) we derive that $f'(k)=0$  if and only if
$2E(k)=F(k)$. The function $F(k)$ is monotonically increasing
as $k \rightarrow 1$, and the function $E(k)$ is monotonically
decreasing at the same time.
We notice now that $2E(k)=F(k)$ for
$k^2 \approx  0.826$ ($E \approx 1.1613, F \approx 2.3181$ for $k^2 =0.825,$
$E \approx 1.1606, F \approx 2.3207$ for $k^2=0.826$, and
$E \approx 1.1599, F \approx 2.3234$ for $k^2=0.827$) (\cite{Tables}).
We can make rigorous estimates but from these tables it is rather evident
that (5.16) holds because $f(\sqrt{0.826}) \approx 0.9352$ and
the estimates given above allow us to make this conclusion
(notice also that $\frac{\pi}{4} < 0.7854$).

Thus we conclude that if potential $p(x)$ which satisfies (6.7) with
$\alpha >0$ induces a torus then the total mean curvature of this torus is
greater than $2\pi^2$.

Thus we conclude that

{\it if potential $p(x)$  induces, via (2.2) and (4.1), mKdV--invariant
Willmore torus of revolution then total squared mean curvature of this torus
is greater or equal to $2\pi^2$ and, moreover, is
equal to $2\pi^2$ only for the Clifford torus.}

The case of mKDV-noninvariant tori will be considered elsewhere.

\bigskip

{\bf Acknowledgement:}
The second author (I.A.T.) was supported by Mathematisches
Forschunginsitut in Oberwolfach and Volkswagenverk Foundation. This work
was done during his stay at the Ruhr-Universitat in Bochum.
He also aknowledges partial support by the Russian Foundation for
Fundamental Studies (grant 94-01-00528).


\newpage

\begin{thebibliography}{99}

\bibitem{BCC}
Bishop A.R., Campbell L.G., and P.G.Channel, Eds.,
{\it Fronts, Interfaces and Patterns}, North--Holland, 1984.

\bibitem{Pelce}
Pelce P., Ed., {\it Dynamics of Curved Fronts}, Academic Press, 1986.

\bibitem{NPW}
Nelson D., Piran T., and Weinberg S., Eds., {\it Statistical
Mechanics of Membranes and Surfaces}, World Scientific, Singapore, 1989.

\bibitem{GPW}
Gross D.G., Piran T., and Weinberg S., Eds., {\it Two--dimensional
Quantum Gravity and Random Surfaces}, World Scientific, Singapore, 1992.

\bibitem{Polyakov}
Polyakov A.M., Phys. Lett. 103  B (1981), 207.

\bibitem{Darboux}
Darboux G., {\it Lecons sur la theorie des surfaces et les applications
geometriques du calcul infinitesimal}, t. 1--4,
Gauthier--Villars, Paris, 1877--1896.

\bibitem{Bianchi}
Bianchi L., {\it Lezioni di geometria differentiale}, 2nd ed.,
Spoerri, Pisa, 1902.

\bibitem{Yau}
Yau S.--T., Ed., {Seminar on differential geometry},
Princeton Univ. Press, Princeton, 1982.

\bibitem{Willmore1982}
Willmore T.J., {\it Total Curvature in Riemannian Geometry},
John Wiley and Sons, New York, 1982.

\bibitem{Willmore1993}
Willmore T.J., {\it Riemannian geometry}, Clarendon Press, Oxford, 1993.

\bibitem{GM}
Ginsparg P., Moore G., {\it Lectures on 2D gravity and 2D string theory},
hep--th 9304011.

\bibitem{ST}
Shiohama K., Takagi A., J. Diff. Geom. 4 (1970), 477--485.

\bibitem{Willmore}
Willmore T.J., J. London Math. Soc. (2) 3 (1971), 307--310.

\bibitem{LY}
Li P. , Yau S.-T., Inventiones Math. 69 (1982), 269--291.

\bibitem{Kusner}
Kusner R., Pacific J. Math. 138 (1989), 317--345.

\bibitem{FPPS}
Ferus D., Pedit F., Pinkall U., Sterling I., J. reine angew. Math.
429 (1992), 1--47.

\bibitem{HJP}
Hertrich--Jeromin U., Pinkall U., J. reine angew. Math.
430 (1992), 21--34.

\bibitem{BB}
Babich M., Bobenko A., Duke Math. J. 72 (1993), 151--185.

\bibitem{Babich}
Babich M., Preprint, Clarkson University, 1994.

\bibitem{Kon}
Konopelchenko B., {\it Induced surfaces and their integrable dynamics},
Studies in Appl. Math. (to appear); preprint BudkerINP 93--114, Novosibirsk,
1993.

\bibitem{Ken}
Kenmotsu K., Math. Ann. 245 (1979), 89--99.

\bibitem{HO}
Hoffmann D.A. and Osserman R., Proceedings of the London Math. Soc.
50 (1985), 27--56.

\bibitem{CK}
Carroll R., Konopelchenko B., {\it Generalized Weierstrass--Enneper
inducing, conformal immersions, and gravity}, to be published.

\bibitem{Bogdan}
Bogdanov L., Teor. Mat. Fyz. 70 (1987), 309.

\bibitem{VN}
Veselov A.P., Novikov S.P., Soviet Math. Dokl. 30 (1984), 588--591,
705--708.

\bibitem{DKN}
Dubrovin B.A., Krichever I.M., Novikov S.P., Soviet Math. Dokl.
17 (1976),  947--951.

\bibitem{Nov}
Zakharov V.E., Manakov S.V., Novikov S.P., and Pitaevskii L.P.,
{\it Theory of solitons}, Nauka, Moscow, 1980.

\bibitem{AS}
Ablowitz M.G., Segur H.,
{\it Solitons and inverse scattering transform}, SIAM, Philadelphia, 1981.

\bibitem{S}
Schief W.K., Nonlinearity 8 (1995), 1--9.

\bibitem{GR}
Gradshtein I.S., Ryzhik I.M., {\it Table of integrals, series and products},
Academic Press, New York--San Francisco--London, 1965.

\bibitem{Tables}
Belyakov V.M., Kravtsova P.I., and Rappoport M.G.,
{\it Tables of elliptic integrals. I.}, Pergamon Press, New York, 1965.

\end{thebibliography}
\end{document}